\documentclass{amsart}
\usepackage[T1]{fontenc}
\usepackage[latin1]{inputenc}
\usepackage{geometry}
\usepackage{setspace}
\usepackage{amssymb}

\makeatletter

\providecommand{\LyX}{L\kern-.1667em\lower.25em\hbox{Y}\kern-.125emX\@}

 \theoremstyle{plain}    
 \newtheorem{thm}{Theorem} 
 \theoremstyle{plain}    
 \newtheorem{cor}{Corollary} 
 \theoremstyle{plain}    
 \newtheorem{lem}{Lemma} 
 \theoremstyle{plain}    
 \newtheorem{prop}{Proposition} 
 \theoremstyle{definition}
 \newtheorem{defn}{Definition}
 \theoremstyle{remark}
 \newtheorem{rem}{Remark}
 \theoremstyle{remark}    
 \newtheorem*{acknowledgement*}{Acknowledgement} 

\usepackage{amsfonts}
\usepackage{amssymb}

\makeatother

\begin{document}

\newcommand{\Q}{\mathbb {Q}}

\newcommand{\Z}{\mathbb {Z}}

\newcommand{\M}[3]{#1 \equiv #2 \, \, \left( \mathrm{mod}\, #3 \right) }

\newcommand{\cyc}[1]{\mathbb {Q}\left[ \zeta _{#1 }\right] }

\newcommand{\fr}[1]{\mathfrak {#1 }}

\newcommand{\id}{\mathbb {I}}

\title{The kernel of the modular representation and the Galois action in RCFT}

\author{P. Bantay}



\curraddr{Institute for Theoretical Physics, Rolland Eotvos University, Budapest}

\email{bantay@poe.elte.hu}

\thanks{Work supported by grant OTKA T19477.}

\date{2001 February 18 }

\begin{abstract}
It is shown that for the modular representations associated to Rational Conformal
Field Theories, the kernel is a congruence subgroup whose level equals the order
of the Dehn-twist. An explicit algebraic characterization of the kernel is given.
It is also shown that the conductor, i.e. the order of the Dehn-twist is bounded
by a function of the number of primary fields, allowing for a systematic enumeration
of the modular representations coming from RCFTs. Restrictions on the spectrum
of the Dehn-twist and arithmetic properties of modular matrix elements are presented.
\end{abstract}
\maketitle

\section{Introduction}

It is known since the work of Cardy \cite{Cardy} that a most important feature
of Conformal Field Theory is modular invariance: the genus one characters of
the chiral algebra afford a unitary representation of the modular group \( \Gamma (1)=SL(2,\mathbb {Z}) \),
and the torus partition function is an invariant sesquilinear combination of
them. In fact, this is not simply a unitary representation, because it comes
with a distinguished basis formed by the genus one characters of the primary
fields. The corresponding representation matrices enjoy many intriguing algebraic
and arithmetic properties, the most important ones being summarized in Verlinde's
theorem \cite{Ver}\cite{MS}, but one may also cite the trace identities of
\cite{TI} and other related results. These special features are related to
the fact that the modular representation provides the defining data of a 3D
Topological Field Theory \cite{TFT}\cite{Frochlich}. 

A most important characteristic of the modular representation is its kernel,
the normal subgroup of \( \Gamma (1) \) whose elements are represented by the
identity matrix. It has been conjectured by several authors \cite{Moore}\cite{Eholzer}\cite{Eh-Sko}\cite{BCIR}
that the kernel is a congruence subgroup, i.e. it contains some principal congruence
subgroup \( \Gamma (N) \). This would have important consequences for a better
understanding of the analytic and arithmetic properties of the genus one characters.
No general proof of the above congruence subgroup property may be found in the
literature, the best known result being that the congruence property holds if
the Dehn-twist \( T \) is represented by a matrix of odd order \cite{CG2},
but for generic RCFTs the order of \( T \) tends to be even, as can be seen
on the example of (Virasoro) minimal models.

The aim of the present paper is to present a proof of the congruence subgroup
property for an arbitrary RCFT. We shall do this by studying the Galois action
\cite{BG}\cite{CG1} in the RCFT and in suitable permutation orbifolds of it.
By exploiting the knowledge of the modular representation of these permutation
orbifolds, we'll be able to give a simple description of the Galois action on
arbitrary modular matrices. This in turn will allow us to write down a simple
expression for the matrix elements of an arbitrary modular transformation, leading
to a straightforward characterization of the kernel and an immediate proof of
the congruence subgroup property. It should be emphasized that the results we
present are much more powerful than the congruence property itself, for they
allow a more refined study of the interplay between the arithmetic properties
of modular matrix elements and the group theoretic properties of the modular
representation.

\section{The Galois action}

In this section we'll review the basic aspects of the Galois action in RCFT,
beginning with some standard facts about the modular representation, mostly
to set the notation. 

Let's consider a Rational Conformal Field Theory \( \mathcal{C} \) with a finite
set \( \mathcal{I} \) of primary fields. The genus one characters \( \chi _{p}\left( \tau \right)  \)
of the primaries \( p\in \mathcal{I} \) transform according to a unitary representation
of the modular group \( \Gamma (1) \), i.e. for \( m=\left( \begin{array}{cc}
a & b\\
c & d
\end{array}\right) \in \Gamma (1) \) we have \begin{equation}
\label{modtrans}
\chi _{p}\left( \frac{a\tau +b}{c\tau +d}\right) =\sum _{q}M_{p}^{q}\chi _{q}\left( \tau \right) 
\end{equation}
In the sequel, we shall always denote by \( M \) the matrix representing \( m\in \Gamma (1) \)
in the basis of genus one characters.

Of special interest are the matrices \( T \) and \( S \) representing \( t=\left( \begin{array}{cc}
1 & 1\\
0 & 1
\end{array}\right)  \) and \( s=\left( \begin{array}{cc}
0 & -1\\
1 & 0
\end{array}\right)  \). As \( s \) and \( t \) generate the modular group \( \Gamma (1) \), any
representation matrix \( M \) may be written in terms of \( S \) and \( T \),
although the resulting expression may be quite cumbersome. It follows from the
defining relations of \( \Gamma (1) \) that \begin{eqnarray}
STS & = & T^{-1}ST^{-1}\label{modrel1} \\
S^{4} & = & 1\label{modrel2} 
\end{eqnarray}
Moreover, it is known that \( S^{2} \) is the charge conjugation operator,
i.e. \begin{equation}
\label{s2}
\left( S^{2}\right) _{p}^{q}=\delta _{p,\overline{q}}
\end{equation}
where \( \overline{q} \) denotes the charge conjugate of the primary \( q \). 

The most important properties of \( S \) and \( T \) are summarized in the
celebrated theorem of Verlinde \cite{Ver}: 

\begin{enumerate}
\item \( T \) is diagonal of finite order.
\item \( S \) is symmetric.
\item The quantities \[
N_{pqr}=\sum _{s\in \mathcal{I}}\frac{S_{ps}S_{qs}S_{rs}}{S_{0s}}\]
are non-negative integers, being the dimension of suitable spaces of holomorphic
blocks. Here and in the sequel, the label \( 0 \) refers to the vacuum of the
theory.
\end{enumerate}
The basic idea in the theory of the Galois action \cite{BG}\cite{CG1} is to
look at the field \( F \) obtained by adjoining to the rationals \( \mathbb {Q} \)
the matrix elements of all modular transformations. One may show that, as a
consequence of Verlinde's theorem, \( F \) is a finite Abelian extension of
\( \mathbb {Q} \). By the celebrated theorem of Kronecker and Weber this means
that \( F \) is a subfield of some cyclotomic field \( \mathbb {Q}\left[ \zeta _{n}\right]  \)
for some integer \( n \), where \( \zeta _{n}=\exp \left( \frac{2\pi i}{n}\right)  \)
is a primitive \( n \)-th root of unity. We'll call the conductor of \( \mathcal{C} \)
the smallest \( n \) for which \( F\subseteq \cyc{n} \) and which is divisible
by the order of the Dehn-twist. 

Among other things, the above results imply that the Galois group \( \mathrm{Gal}\left( F/\mathbb {Q}\right)  \)
is a homomorphic image of the Galois group \( \mathcal{G}_{n}=\mathrm{Gal}\left( \cyc{n}/\mathbb {Q}\right)  \).
But it is well known that \( \mathcal{G}_{n} \) is isomorphic to the group
\( \left( \Z /n\Z \right) ^{*} \) of prime residues modulo \( n \), its elements
being the Frobenius maps \( \sigma _{l}:\cyc{n}\rightarrow \cyc{n} \) that
leave \( \Q  \) fixed, and send \( \zeta _{n} \) to \( \zeta _{n}^{l} \)
for \( l \) coprime to \( n \). Consequently, the maps \( \sigma _{l} \)
are automorphisms of \( F \) over \( \Q  \).

According to \cite{CG1}, we have (for \( l \) coprime to the conductor) \begin{equation}
\label{sls}
\sigma _{l}\left( S_{p}^{q}\right) =\varepsilon _{l}(q)S_{p}^{\pi _{l}q}
\end{equation}
for some permutation \( \pi _{l}\in Sym\left( \mathcal{I}\right)  \) of the
primaries and some function \( \varepsilon _{l}:\mathcal{I}\rightarrow \left\{ -1,+1\right\}  \).
In other words, upon introducing the orthogonal monomial matrices \begin{equation}
\label{gldef}
\left( G_{l}\right) _{p}^{q}=\varepsilon _{l}(q)\delta _{p,\pi _{l}q}
\end{equation}
 and denoting by \( \sigma _{l}\left( M\right)  \) the matrix that one obtains
by applying \( \sigma _{l} \) to \( M \) elementwise, we have \begin{equation}
\label{slg}
\sigma _{l}\left( S\right) =SG_{l}=G_{l}^{-1}S
\end{equation}
Note that for \( l \) and \( m \) both coprime to the conductor \begin{eqnarray*}
\pi _{lm} & = & \pi _{l}\pi _{m}\\
G_{lm} & = & G_{l}G_{m}
\end{eqnarray*}
 The Galois action on \( T \) is even simpler, for \( T \) is diagonal, and
its eigenvalues are roots of unity, consequently \begin{equation}
\label{slt}
\sigma _{l}\left( T\right) =T^{l}
\end{equation}

\section{Galois action on \protect\( \Lambda \protect \) matrices}

In this section we'll study the Galois action in some appropriate permutation
orbifold \cite{KS}. According to the Orbifold Covariance Principle \cite{OCP},
all the properties of the Galois action reviewed in the previous section should
hold in the permutation orbifold, in particular the Galois action on the \( S \)-matrix
elements may be described via suitable permutations \( \pi _{l} \) of the primaries
of the orbifold and signs \( \varepsilon _{l} \) . This will in turn allow
us to determine the Galois action on \( \Lambda  \)-matrices.

We fix a positive integer \( N \), and consider the group \( \Omega  \) generated
by the cyclic permutation \( \left( 1,\ldots ,N\right)  \). One can then form
the permutation orbifold \( \mathcal{C}\wr \Omega  \) according to \cite{BHS}\cite{PO1},
which is a new RCFT with explicitly known genus one characters and modular transformation
matrices. Among the primary fields of the permutation orbifold \( \mathcal{C}\wr \Omega  \)
there is a subset \( \mathcal{J} \) of special relevance to us. The primaries
in \( \mathcal{J} \) are labeled by triples \( \left[ p,n,k\right]  \), where
\( p\in \mathcal{I} \) is a primary of \( \mathcal{C} \), while \( n \) and
\( k \) are integers mod \( N \). The subset of those \( \left[ p,n,k\right]  \)
where \( n \) is coprime to \( N \) will be denoted by \( \mathcal{J}_{0} \).
It follows from the general theory of permutation orbifolds \cite{BHS}\cite{PO1}
that the primary fields \( \left[ p,n,k\right] \in \mathcal{J}_{0} \) have
vanishing \( S \)-matrix elements with the primaries not in \( \mathcal{J} \),
while for \( \left[ q,m,l\right] \in \mathcal{J} \) we have\footnote{%
For the definition and basic properties of \( \Lambda  \)-matrices, see the
Appendix.
}\begin{equation}
\label{Smat}
S_{\left[ p,n,k\right] }^{\left[ q,m,l\right] }=\frac{1}{N}\zeta _{N}^{-(km+ln)}\Lambda _{p}^{q}\left( \frac{m\hat{n}}{N}\right) 
\end{equation}
where \( \hat{n} \) denotes the mod \( N \) inverse of \( n \) and \( \zeta _{N}=\exp \left( \frac{2\pi i}{N}\right)  \). 

Let's now consider the action of a Galois transformation \( \sigma _{l} \)
with \( l \) coprime to \( N \). According to the results of the previous
section, we have \footnote{%
We'll always assume that \( l \) is coprime to the conductor, but this is no
loss of generality, as it can always be achieved according to Dirichlet's theorem
on primes in arithmetic progressions.
} \begin{equation}
\label{gal1}
\sigma _{l}\left( S_{q}^{\left[ p,n,k\right] }\right) =\tilde{\varepsilon }_{l}\left( p,n,k\right) S_{q}^{\tilde{\pi }_{l}\left[ p,n,k\right] }
\end{equation}
for some permutation \( \tilde{\pi }_{l} \) of the primaries of \( \mathcal{C}\wr \Omega  \)
and some signs \( \tilde{\varepsilon }_{l} \). 

\begin{lem}
The set \( \mathcal{J} \) is invariant under the permutations \( \tilde{\pi }_{l} \),
i.e. \( \tilde{\pi }_{l}\left( \mathcal{J}\right) =\mathcal{J} \). For a primary
\( \left[ p,n,k\right] \in \mathcal{J}_{0} \) one has \begin{equation}
\label{pi}
\tilde{\pi }_{l}\left[ p,n,k\right] =\left[ \pi _{l}p,ln,\tilde{k}\right] 
\end{equation}
for some function \( \tilde{k} \) of \( l,p,n \) and \( k \), and \begin{equation}
\label{eps}
\tilde{\varepsilon }_{l}\left( p,n,k\right) =\varepsilon _{l}(p)
\end{equation}

\end{lem}
\begin{proof}
First, let's fix \( \left[ p,n,k\right] \in \mathcal{J} \). According to Eq.(\ref{Smat}),
we have \[
S^{\left[ p,n,k\right] }_{\left[ q,1,0\right] }=\frac{1}{N}\zeta _{N}^{-k}\Lambda _{q}^{p}\left( \frac{n}{N}\right) \]
and this expression differs from \( 0 \) for at least one \( q\in \mathcal{I} \),
by the unitarity of \( \Lambda  \)-matrices. Select such a \( q\in \mathcal{I} \),
and apply \( \sigma _{l} \) to both sides of the equation. One gets that \[
\tilde{\varepsilon }_{l}\left( p,n,k\right) S_{\left[ q,1,0\right] }^{\tilde{\pi }_{l}\left[ p,n,k\right] }=\sigma _{l}\left( S_{\left[ q,1,0\right] }^{\left[ p,n,k\right] }\right) \]
differs from \( 0 \), but this can only happen if \( \tilde{\pi }_{l}\left[ p,n,k\right] \in \mathcal{J} \)
because \( \left[ q,1,0\right] \in \mathcal{J}_{0} \).

Next, for \( \left[ p,n,k\right] \in \mathcal{J}_{0} \) consider\begin{equation}
\label{S0}
S_{\left[ q,0,m\right] }^{\left[ p,n,k\right] }=\frac{1}{N}\zeta _{N}^{-nm}S_{q}^{p}
\end{equation}
Applying \( \sigma _{l} \) to both sides of the above equation we get from
Eq.(\ref{sls}) \[
\tilde{\varepsilon }_{l}(p,n,k)S_{\left[ q,0,m\right] }^{\tilde{\pi }_{l}\left[ p,n,k\right] }=\frac{1}{N}\zeta _{N}^{-lnm}\varepsilon _{l}(p)S_{q}^{\pi _{l}p}\]
But the lhs. equals \[
\tilde{\varepsilon }_{l}(p,n,k)\frac{1}{N}\zeta _{N}^{-\tilde{n}m}S_{q}^{\tilde{p}}\]
 according to Eq.(\ref{S0}) if \( \tilde{\pi }_{l}\left[ p,n,k\right] =\left[ \tilde{p},\tilde{n},\tilde{k}\right]  \).
Equating both sides we arrive at \[
S_{q}^{\tilde{p}}=\varepsilon _{l}(p)\tilde{\varepsilon }_{l}(p,n,k)\zeta _{N}^{-m(\tilde{n}-ln)}S_{q}^{\pi _{l}p}\]
 As the lhs. is independent of \( m \), we must have \[
\M{\tilde{n}}{ln}{N}\]
and \[
\tilde{p}=\pi _{l}p\]
 as well as \[
\tilde{\varepsilon }_{l}(p,n,k)=\varepsilon _{l}(p)\]
 
\end{proof}
\begin{prop}
For \( l \) coprime to the denominator of \( r \) \begin{equation}
\label{Lmat2}
\sigma _{l}\left( \Lambda \left( r\right) \right) =\Lambda \left( lr\right) G_{l}Z_{l}(r^{*})=Z_{l}\left( r\right) G_{l}^{-1}\Lambda (\hat{l}r)
\end{equation}
where \( Z_{l}(r) \) is a diagonal matrix whose order divides the denominator
\( N \) of \( r \), and \( \hat{l} \) is the mod \( N \) inverse of \( l \).
\end{prop}
\begin{proof}
Write \( r=\frac{n}{N} \), with \( n \) coprime to \( N \), and recall that
\( r^{*}=\frac{\hat{n}}{N} \), where \( \hat{n} \) is the mod \( N \) inverse
of \( n \). Consider \begin{equation}
\label{S1}
S_{\left[ q,1,m\right] }^{\left[ p,n,k\right] }=\frac{1}{N}\zeta _{N}^{-(k+nm)}\Lambda _{q}^{p}\left( \frac{n}{N}\right) 
\end{equation}
Applying \( \sigma _{l} \) to both sides of the above equation and taking into
account Eqs. (\ref{pi},\ref{eps}), we get \[
\varepsilon _{l}(p)S_{\left[ q,1,m\right] }^{\left[ \pi _{l}p,ln,\tilde{k}\right] }=\frac{1}{N}\zeta _{N}^{-l(k+nm)}\sigma _{l}\left( \Lambda _{q}^{p}\left( \frac{n}{N}\right) \right) \]
 But form Eq.(\ref{S1}), the lhs. equals \[
\varepsilon _{l}(p)\frac{1}{N}\zeta _{N}^{-(\tilde{k}+lnm)}\Lambda _{q}^{\pi _{l}p}\left( \frac{ln}{N}\right) \]
so after rearranging, we get \begin{equation}
\label{Lmat1}
\sigma _{l}\left( \Lambda _{q}^{p}\left( \frac{n}{N}\right) \right) =\varepsilon _{l}(p)\zeta _{N}^{(lk-\tilde{k})}\Lambda _{q}^{\pi _{l}p}\left( \frac{ln}{N}\right) 
\end{equation}
As the lhs. is independent of \( k \) \begin{equation}
\label{k_0}
\M{\tilde{k}-lk}{k_{0}\left( l,p,r\right) }{N}
\end{equation}
Of course, the point is that \( k_{0} \) is independent of \( k \). If we
introduce the diagonal matrix \begin{equation}
\label{Zmat}
Z_{l}\left( r\right) _{p}^{q}=\delta _{p,q}\zeta _{N}^{-k_{0}(l,p,r^{*})}
\end{equation}
then it is obvious that the order of \( Z_{l}(r) \) divides the denominator
of \( r \) and the first equality in Eq.(\ref{Lmat2}) holds. The second follows
simply from \( \Lambda \left( r\right) ^{T}=\Lambda \left( r^{*}\right)  \).
\end{proof}
Note that it follows from the definition of \( Z_{l}(r) \) that \( Z_{l}(0)=\mathbb {I} \)
and \( Z_{l}(r+1)=Z_{l}(r) \).

\begin{lem}
\begin{equation}
\label{zcoc}
G_{l}^{-1}Z_{m}(\hat{l}r)G_{l}=Z_{lm}(r)Z_{l}^{-m}(r)
\end{equation}
whenever both \( l \) and \( m \) are coprime to the denominator of \( r \).
\end{lem}
\begin{proof}
This follows at once from \( \sigma _{lm}=\sigma _{l}\sigma _{m} \) applied
to \( \Lambda (r) \).
\end{proof}
\begin{lem}
If \( n \) is coprime to the denominator of \( r \), then \begin{equation}
\label{zmult}
Z_{l}^{n}\left( r\right) =Z_{l}\left( nr\right) 
\end{equation}

\end{lem}
\begin{proof}
Let \( r=\frac{m}{N} \), with \( m \) coprime to \( N \). According to Eq.(\ref{Smat})
\[
NS_{\left[ p,n,0\right] }^{\left[ q,nm,0\right] }=\Lambda _{p}^{q}\left( r\right) \]
Applying \( \sigma _{l} \) to both sides we get \[
\sigma _{l}\left( \Lambda _{p}^{q}\left( r\right) \right) =N\varepsilon _{l}(q)S_{\left[ p,n,0\right] }^{\left[ \pi _{l}q,lnm,k_{0}\right] }=\varepsilon _{l}(q)\zeta _{N}^{-nk_{0}\left( l,q,\frac{nm}{N}\right) }\Lambda \left( lr\right) _{p}^{\pi _{l}q}\]
i.e. \[
\Lambda \left( lr\right) G_{l}Z_{l}\left( r^{*}\right) =\Lambda \left( lr\right) G_{l}Z_{l}^{n}\left( \hat{n}r^{*}\right) \]
by Eq.(\ref{Lmat2}), i.e. \[
Z_{l}\left( r\right) =Z_{l}^{n}\left( \hat{n}r\right) \]
But this is equivalent to the assertion. 
\end{proof}
\begin{thm}
For all \( l \) coprime to the conductor\begin{equation}
\label{gtcom}
G_{l}^{-1}TG_{l}=T^{l^{2}}
\end{equation}

\end{thm}
\begin{proof}
Take \( N \) equal to the order of \( T \). Clearly, \( N \) divides the
conductor. By Eq.(\ref{Lambda1n}) of the Appendix, \[
\Lambda \left( \frac{1}{N}\right) =T^{-\frac{1}{N}}S^{-1}T^{-N}ST^{-\frac{1}{N}}=T^{-\frac{2}{N}}\]
 But \[
\sigma _{l}^{2}\left( \Lambda \left( \frac{1}{N}\right) \right) =\sigma _{l}\left( \Lambda \left( \frac{l}{N}\right) G_{l}Z_{l}\left( \frac{1}{N}\right) \right) =Z_{l}\left( \frac{l}{N}\right) G_{l}^{-1}\Lambda \left( \frac{1}{N}\right) G_{l}Z_{l}^{l}\left( \frac{1}{N}\right) \]
 or in other words \[
T^{-\frac{2l^{2}}{N}}=G_{l}^{-1}T^{-\frac{2}{N}}G_{l}Z_{l}^{2}\left( \frac{l}{N}\right) \]
because both \( Z_{l} \) and \( T \) are diagonal, hence they commute. If
\( N \) is even, then taking the \( \frac{N}{2} \)-th power of the last equation
gives the result. If \( N \) is odd, then there exists \( k \) such that \( \M{2k}{1}{N} \),
and taking the \( kN \)-th power of the equation gives again Eq.(\ref{gtcom}). 
\end{proof}
\begin{rem}
This result has been conjectured in \cite{CG2}, where some of its consequences
(to be reviewed in the next section) had been derived.
\end{rem}
\begin{cor}
\begin{equation}
\label{gcom}
G_{l}^{-1}MG_{l}=\sigma _{l}^{2}\left( M\right) 
\end{equation}
\label{gcomlemma}for any modular matrix \( M \). 
\end{cor}
\begin{proof}
We have just seen that it holds for \( T \), and it also holds for \( S \)
by Eq.(\ref{slg}). But \( s \) and \( t \) generate \( \Gamma (1) \), consequently
the claim should hold for any \( m\in \Gamma (1) \).
\end{proof}
\begin{prop}
If \( l \) is coprime to both the conductor and the denominator \label{gtcom1lemma}of
\( r \), then \begin{equation}
\label{gtcom1}
G_{l}^{-1}T^{r}G_{l}=T^{l^{2}r}Z_{l}^{l}\left( r\right) 
\end{equation}

\end{prop}
\begin{proof}
Let's write \( r=\frac{n}{N} \). We know from \cite{BHS}\cite{PO1} that \[
\omega _{\left[ p,n,k\right] }=\zeta _{N}^{nk}\omega _{p}^{1/N}\]
where \( \omega _{p}=\exp \left( 2\pi i\left( \Delta _{p}-\frac{c}{24}\right) \right)  \)
is the exponentiated conformal weight of the primary \( p \), i.e. the matrix
element \( T_{pp} \). By Eq.(\ref{gtcom}) \[
\omega _{p}^{l^{2}/N}=\omega _{\tilde{\pi }_{l}\left[ p,1,0\right] }=\omega _{\left[ \pi _{l}p,l,k_{0}\right] }=\zeta _{N}^{lk_{0}}\omega _{\pi _{l}p}^{1/N}\]
 i.e. \[
Z^{l}_{l}\left( \frac{1}{N}\right) _{pp}=\zeta _{N}^{-lk_{0}(l,p,1/N)}=\omega _{\pi _{l}p}^{1/N}\omega _{p}^{-l^{2}/N}\]
 or in other words \[
Z^{l}_{l}\left( \frac{1}{N}\right) =G_{l}^{-1}T^{\frac{1}{N}}G_{l}T^{-\frac{l^{2}}{N}}\]
Taking the \( n \)-th power of the last equation and using Eq.(\ref{zmult})
gives the result, because \( Z_{l} \) and \( T \) commute. 
\end{proof}
\begin{lem}
Suppose that \( l \) is coprime to the denominator of both \( r_{1} \) and
\( r_{2} \). Then\begin{equation}
\label{zadd}
Z_{l}(r_{1})Z_{l}(r_{2})=Z_{l}\left( r_{1}+r_{2}\right) 
\end{equation}
 
\end{lem}
\begin{proof}
\[
T^{l^{2}(r_{1}+r_{2})}Z_{l}^{l}(r_{1})Z_{l}^{l}(r_{2})=G_{l}^{-1}T^{r_{1}+r_{2}}G_{l}=T^{l^{2}(r_{1}+r_{2})}Z_{l}^{l}\left( r_{1}+r_{2}\right) \]
according to Proposition \ref{gtcom1lemma}.
\end{proof}

\section{The Galois action revisited}

In this section we translate the results of the previous section about the Galois
action on \( \Lambda  \)-matrices to statements about the Galois action on
modular matrix elements. 

\begin{prop}
For \( l \) coprime to the conductor, \begin{equation}
\label{Gl}
G_{l}=S^{-1}T^{l}ST^{\hat{l}}ST^{l}
\end{equation}
where \( \hat{l} \) denotes the inverse of \( l \) modulo the conductor.
\end{prop}
\begin{proof}
Let's apply \( \sigma _{l} \) to both sides of the modular relation Eq.(\ref{modrel1}).
One gets\[
SG_{l}T^{l}G_{l}^{-1}S=T^{-l}SG_{l}T^{-l}\]
After rearranging and using Eq.(\ref{gtcom}) we get the assertion. 
\end{proof}
\begin{prop}
\label{condlemma}The conductor equals the order \( N \) of \( T \), and \( F=\cyc{N} \).
\end{prop}
\begin{proof}
Clearly, \( N \) divides the conductor, for the eigenvalues of \( T \) are
roots of unity. But for \( \M{l}{1}{N} \) we have \[
\sigma _{l}\left( T\right) =T^{l}=T\]
and \[
\sigma _{l}\left( S\right) =T^{l}ST^{\hat{l}}ST^{l}=S\]
according to Eqs.(\ref{slt},\ref{slg},\ref{Gl}), consequently all such \( l \)
fixes \( F \), i.e. \( F\subseteq \cyc{N} \). On the other hand, if \( \sigma _{l} \)
fixes \( F \) then \( \M{l}{1}{N} \), i.e. \( Gal\left( F/\Q \right) =Gal\left( \cyc{N}/\Q \right)  \). 
\end{proof}
\begin{rem}
The above two propositions had been derived in\cite{CG2}, upon postulating
Eq.(\ref{gtcom}). Earlier, they have been conjectured in \cite{Bauer}.
\end{rem}
\begin{lem}
\label{Gdiag}If \( G_{l} \) is diagonal, then \( G_{l}=\pm \mathbb {I} \).
\end{lem}
\begin{proof}
\( G_{l} \) diagonal means that \( \pi _{l} \) is the identity permutation.
If we apply \( \sigma _{l} \) to \( S_{0p} \), we get \[
\sigma _{l}\left( S_{0p}\right) =\varepsilon _{l}\left( 0\right) S_{0p}=\varepsilon _{l}\left( p\right) S_{0p}\]
because \( \pi _{l}p=p \). But \( S_{0p}>0 \), consequently \( \varepsilon _{l}(p)=\varepsilon _{l}(0)=\pm 1 \),
proving the lemma. 
\end{proof}
\begin{prop}
Let \( N_{0} \) denote the order of the matrix \( \omega _{0}^{-1}T \), i.e.
the least common multiple of the denominators of the conformal weights. Then
\( N=eN_{0} \), where the integer \( e \) divides \( 12 \). Moreover, the
greatest common divisor of \( e \) and \( N_{0} \) is either 1 or 2.
\end{prop}
\begin{proof}
That \( N_{0} \) divides \( N \) is obvious. If we consider the matrix \( T^{N_{0}} \),
then it will equal \( \zeta =\omega _{0}^{N_{0}} \) times the identity matrix.
But \( N \) is the exact order of \( T \), consequently \( \zeta  \) is a
primitive \( e \)-th root of unity. On the other hand, Eq.(\ref{gtcom}) implies
\[
G_{l}^{-1}T^{N_{0}}G_{l}=T^{l^{2}N_{0}}\]
for all \( l \) coprime to \( N \), i.e. \[
\zeta ^{l^{2}}=\zeta \]
In other words, the exponent of \( \left( \Z /e\Z \right) ^{*} \) should divide
\( 2 \), and this can only happen if \( e \) divides \( 24 \).

Next, let's consider \( l=1+6N_{0} \). By the above result, \( l \) is coprime
to both \( N_{0} \) and \( e \), consequently it is also coprime to \( N \).
If \( \hat{l} \) denotes its inverse mod \( N \), then we have \[
\hat{l}=\left\{ \begin{array}{cc}
1+6N_{0} & \textrm{if }\, \, N_{0}\, \, \textrm{ is odd},\\
1-6N_{0} & \textrm{if }\, \, N_{0}\, \, \textrm{ is even}.
\end{array}\right. \]
 According to Eq.(\ref{Gl}), \[
G_{l}=ST^{l}S^{-1}T^{\hat{l}}ST^{l}=\zeta ^{6\left( 2\pm 1\right) }\id \]
But if \( G_{l} \) is diagonal, then it equals \( \pm \id  \), so that \( \zeta ^{12\left( 2\pm 1\right) }=1 \).
Together with \( \zeta ^{24}=1 \), this yields \( \zeta ^{12}=1 \).

Finally, let \( D \) denote the gcd of \( e \) and \( N_{0} \), and consider
\( l=1+\frac{N}{D} \). It is straightforward that \( l \) is coprime to \( N \),
and its inverse modulo \( N \) is \( \hat{l}=1-\frac{N}{D} \). But \[
G_{l}=\zeta ^{\frac{e}{D}}\id \]
 and once again this should equal \( \pm \id  \), consequently \( D \) divides
2.
\end{proof}
\begin{rem}
As it turns out, there are no further restrictions on the value of \( N/N_{0} \),
i.e. each divisor of \( 12 \) occurs for some RCFT. But for a given value of
\( N/N_{0} \), techniques similar to the above yield more restrictions on \( N_{0} \),
e.g. \( N/N_{0}=12 \) cannot occur unless \( \M{N_{0}}{\pm 1}{6} \). 
\end{rem}
\begin{cor}
\( N_{0} \) times the central charge is an even integer.
\end{cor}
\begin{proof}
If \( c \) denotes the central charge, then \( \omega _{0}=\exp \left( -\pi i\frac{c}{12}\right)  \).
But \( \omega _{0}^{N}=1 \), so by the above \( \omega _{0}^{12N_{0}}=1 \),
proving the claim.
\end{proof}
\begin{prop}
\label{condbound}There exists a function \( N(r) \) such that the conductor
\( N \) divides \( N(r) \) if the number of primary fields - i.e. the dimension
of the modular representation - is \( r \). 
\end{prop}
\begin{proof}
Let's denote by \( s(r) \) the exponent of the symmetric group \( S_{r} \)
of degree \( r \). As the permutations \( \pi _{l} \) belong to \( S_{r} \),
it follows that \( \pi _{l}^{s(r)} \) is the identity permutation, i.e. \( G_{l}^{s(r)} \)
is diagonal. By Lemma \ref{Gdiag} this means that \( G_{l}^{s(r)}=\pm \mathbb {I} \),
which in turn implies that \( \M{l^{2s(r)}}{1}{N} \) for all \( l \), i.e.
the exponent of the group \( \left( \Z /N\Z \right) ^{*} \) divides \( 2s(r) \).
Taking \( N(r) \) to be the greatest integer satisfying this last condition
proves the proposition.
\end{proof}
According to this proposition, for a given number of primary fields one has
only a finite number of consistent choices for the conductor \( N \) and the
matrix \( T \). The upper bound for \( N(r) \) given in the proof may be greatly
improved by exploiting the known properties of the matrices \( G_{l} \) (e.g.
that they commute), which leads to the following table for small values of \( r \)
: 

\vspace{0.375cm}
{\centering \begin{tabular}{|c|c|}
\hline 
\( r \) &
\( N(r) \) \\
\hline 
\hline 
2&
240\\
\hline 
3&
5040\\
\hline 
4&
10080\\
\hline 
5&
1441440\\
\hline 
\end{tabular}\par}
\vspace{0.375cm}

\section{The kernel of the modular representation}

Recall that the kernel \( \mathcal{K} \) consists of those modular transformations
which are represented by the identity matrix, i.e. \[
\mathcal{K}=\left\{ m\in \Gamma (1)\, |\, M_{p}^{q}=\delta _{p,q}\right\} \]

\begin{prop}
Let \( m=\left( \begin{array}{cc}
a & b\\
c & d
\end{array}\right) \in SL(2,\mathbb {Z}) \), \( l \) coprime to the conductor, \( lk=1+\alpha c \) and \[
\hat{m}=\left( \begin{array}{cc}
la & b+\alpha ad\\
c & kd
\end{array}\right) \]
 Then \begin{equation}
\label{modrep0}
\sigma _{l}\left( M\right) =\hat{M}G_{l}T^{-\alpha ld}
\end{equation}

\end{prop}
\begin{proof}
According to Eq.(\ref{Lmat2}) \[
\sigma _{l}\left( M\right) =\sigma _{l}\left( T^{a/c}\Lambda \left( \frac{a}{c}\right) T^{d/c}\right) =T^{al/c}\Lambda \left( \frac{al}{c}\right) G_{l}Z_{l}\left( \frac{d}{c}\right) T^{dl/c}\]
because \( \left( \frac{a}{c}\right) ^{*}=\frac{d}{c} \). But 

\[
\Lambda \left( \frac{al}{c}\right) =T^{-al/c}\hat{M}T^{-kd/c}\]
so \[
\sigma _{l}\left( M\right) =\hat{M}T^{-kd/c}G_{l}Z_{l}\left( \frac{d}{c}\right) T^{ld/c}\]
From Eq.(\ref{gtcom1}) \[
T^{-kd/c}G_{l}=G_{l}T^{-l^{2}kd/c}Z_{l}\left( \frac{-lkd}{c}\right) \]
Putting all this together, we get \[
\sigma _{l}\left( M\right) =\hat{M}G_{l}T^{-ld\left( 1+\alpha c\right) /c}Z_{l}\left( \frac{-d}{c}\right) Z_{l}\left( \frac{d}{c}\right) T^{ld/c}\]
which proves the claim according to Eq.(\ref{zadd}). 
\end{proof}
\begin{cor}
If \( m=\left( \begin{array}{cc}
a & b\\
c & d
\end{array}\right) \in \Gamma (1) \) with \( d \) coprime to the conductor, then

\begin{equation}
\label{modrep1}
\sigma _{d}\left( M\right) =T^{b}S^{-1}T^{-c}\sigma _{d}\left( S\right) 
\end{equation}
 
\end{cor}
\begin{proof}
If we take \( l=d \) and \( k=a \) in Eq.(\ref{modrep0}), then \( \alpha =b \)
and \( \hat{m}=t^{b}s^{-1}t^{-c}st^{b} \).
\end{proof}
\begin{thm}
Let \( d \) be coprime to the conductor \( N \). Then \( \left( \begin{array}{cc}
a & b\\
c & d
\end{array}\right) \in \Gamma (1) \) belongs to the kernel \( \mathcal{K} \) if and only if \begin{equation}
\label{kern1}
\sigma _{d}\left( S\right) T^{b}=T^{c}S
\end{equation}

\end{thm}
\begin{proof}
If \( m\in \mathcal{K} \), then \( \sigma _{d}\left( M\right)  \) is the identity
matrix, and Eq.(\ref{modrep1}) implies the result. 
\end{proof}
\begin{rem}
Apparently, the above criterion does only apply in case \( d \) is coprime
to the conductor. But if \( m\in \mathcal{K} \), then \[
t^{-k}mt^{k}=\left( \begin{array}{cc}
a-kc & b-k^{2}c+k(a-d)\\
c & d+kc
\end{array}\right) \]
also belongs to \( \mathcal{K} \). Dirichlet's theorem on primes in arithmetic
progressions ensures that \( d+kc \) is coprime to the conductor for infinitely
many \( k \), reducing the general case to the coprime one.
\end{rem}
For the next result, recall that \[
\Gamma _{1}\left( N\right) =\left\{ \left( \begin{array}{cc}
a & b\\
c & d
\end{array}\right) \in \Gamma \left( 1\right) \, |\, \M{a,d}{1}{N},\, \, \M{c}{0}{N}\right\} \]
and\[
\Gamma \left( N\right) =\left\{ \left( \begin{array}{cc}
a & b\\
c & d
\end{array}\right) \in \Gamma _{1}\left( N\right) \, |\, \M{b}{0}{N}\right\} \]

\begin{thm}
\[
\mathcal{K}\cap \Gamma _{1}\left( N\right) =\Gamma \left( N\right) \]
 In particular, \( \mathcal{K} \) is a congruence subgroup of level \( N \). 
\end{thm}
\begin{proof}
If \( \left( \begin{array}{cc}
a & b\\
c & d
\end{array}\right) \in \Gamma _{1}\left( N\right)  \), then Eq.(\ref{kern1}) reduces to \( ST^{b}=S \), whose only solution is
\( \M{b}{0}{N} \).
\end{proof}
The above result means that the modular representation factors through \( SL_{2}\left( N\right) \cong \Gamma (1)/\Gamma (N) \),
i.e. there exists a representation \( D \) of \( SL_{2}\left( N\right)  \)
such that the modular representation is the composite map \( D\circ \mu _{N} \),
where we denote by \( \mu _{N} \) the natural homomorphism \( \Gamma (1)\rightarrow SL_{2}\left( N\right)  \).
The representation \( D \) gives us an elegant description of the Galois action.

\begin{defn}
For \( l \) coprime to \( N \), define the automorphism \( \tau _{l}:SL_{2}\left( N\right) \rightarrow SL_{2}\left( N\right)  \)
by \begin{equation}
\label{taudef}
\tau _{l}\left( \begin{array}{cc}
a & b\\
c & d
\end{array}\right) =\left( \begin{array}{cc}
a & lb\\
\hat{l}c & d
\end{array}\right) 
\end{equation}
where \( \hat{l} \) is the mod \( N \) inverse of \( l \).
\end{defn}
\begin{thm}
\begin{equation}
\label{gal2}
\sigma _{l}\circ D=D\circ \tau _{l}
\end{equation}

\end{thm}
\begin{proof}
As a homomorphic image of \( \Gamma (1) \), \( SL_{2}\left( N\right)  \) is
generated by \( \mu _{N}(s) \) and \( \mu _{N}\left( t\right)  \), so it is
enough to verify Eq.(\ref{gal2}) for these matrices. For \( t \), the lhs.
reads \[
\sigma _{l}\circ D\circ \mu _{N}(t)=T^{l}\]
while the rhs. is \[
D\left( \begin{array}{cc}
1 & l\\
0 & 1
\end{array}\right) =T^{l}\]
and these are clearly equal. For the case of \( s \), the lhs. reads \[
\sigma _{l}\left( S\right) =SG_{l}=T^{l}ST^{\hat{l}}ST^{l}\]
according to Eq.(\ref{Gl}), while the rhs. is 

\[
D\left( \begin{array}{cc}
0 & -l\\
\hat{l} & 0
\end{array}\right) \]
 It is easy to check that \[
\M{\left( \begin{array}{cc}
0 & -l\\
\hat{l} & 0
\end{array}\right) }{t^{l}st^{\hat{l}}st^{l}}{N}\]
 proving the claim. 
\end{proof}
Note that in terms of the representation \( D \) we have \[
G_{l}=D\left( \begin{array}{cc}
\hat{l} & 0\\
0 & l
\end{array}\right) \]
i.e. the matrices \( G_{l} \) represent the diagonal subgroup of \( SL_{2}\left( N\right)  \).

Without going into the details, we note without proof that Eq.(\ref{kern1})
implies the following results :

\begin{enumerate}
\item \( G_{d}^{2}=1 \), consequently \( \M{d^{4}}{1}{N} \).
\item \( \M{b,c}{1-d^{2}}{N_{0}} \) .
\item \( \varepsilon _{d}\left( 0\right) =\omega _{0}^{c-b} \), in particular \( \M{2c}{2b}{N} \).
\item \( \M{2c}{0}{N_{0}} \).
\end{enumerate}

\section{Examples}

Let's look at a couple of simple examples to illustrate the above results. First,
let's consider the Lee-Yang model, i.e. the minimal model \( \mathcal{M}\left( 5,2\right)  \).
It has central charge \( c=-\frac{22}{5} \) and two primary fields, whose conformal
weights are \( 0 \) and \( -\frac{1}{5} \). The \( S \) matrix reads \[
S=\frac{2}{\sqrt{5}}\left( \begin{array}{cc}
-\sin \left( \frac{2\pi }{5}\right)  & \sin \left( \frac{4\pi }{5}\right) \\
\sin \left( \frac{4\pi }{5}\right)  & \sin \left( \frac{2\pi }{5}\right) 
\end{array}\right) \]

The conductor \( N \) equals \( 60 \), while \( N_{0}=5 \), providing an
example where the ratio \( N/N_{0} \) attains the upper bound of Lemma. As
it turns out, the image \( \overline{\mathcal{K}}=\mu _{N}(\mathcal{K}) \)
of the kernel is a non-Abelian group of order \( 192 \), whose center is \( \Z _{2}^{2} \)
, and whose derived subgroup is \( \Z _{2}^{3} \). A small generating set for
\( \overline{\mathcal{K}} \) is provided by the matrices \[
\left( \begin{array}{cc}
19 & 5\\
5 & 14
\end{array}\right) ,\left( \begin{array}{cc}
31 & 35\\
5 & 56
\end{array}\right) ,\left( \begin{array}{cc}
56 & 5\\
35 & 31
\end{array}\right) \]

As a second example, let's consider the Ising model, i.e. the minimal model
\( \mathcal{M}\left( 4,3\right)  \). The central charge is \( c=\frac{1}{2} \),
there are \( 3 \) primary fields of conformal weights \( 0,\frac{1}{2} \)
and \( \frac{1}{16} \), and \[
S=\frac{1}{2}\left( \begin{array}{ccc}
1 & 1 & \sqrt{2}\\
1 & 1 & -\sqrt{2}\\
\sqrt{2} & -\sqrt{2} & 0
\end{array}\right) \]
The conductor equals \( 48 \), while \( N_{0}=16 \). In this case \( \overline{\mathcal{K}} \)
has order \( 64 \), it's derived subgroup is \( \Z _{2} \), while its center
is \( \Z _{2}^{2}\times \Z _{4} \). A small generating set for \( \overline{\mathcal{K}} \)
is given by the \( SL_{2}\left( 48\right)  \) matrices \[
\left( \begin{array}{cc}
43 & 40\\
40 & 35
\end{array}\right) ,\left( \begin{array}{cc}
29 & 40\\
40 & 37
\end{array}\right) ,\left( \begin{array}{cc}
21 & 8\\
40 & 45
\end{array}\right) ,\left( \begin{array}{cc}
35 & 40\\
40 & 43
\end{array}\right) \]

Finally, the following table contains the relevant data for some minimal models,
the \( pq \)th entry giving respectively the conductor \( N \), the ratio
\( N/N_{0} \), and the index \( \left[ \mathcal{K}:\Gamma \left( N\right) \right] =\left| \mu _{N}\left( \mathcal{K}\right) \right|  \)
for the minimal model \( \mathcal{M}(q,p) \) (the first entry in the first
column being the Lee-Yang model discussed above). 

\vspace{0.3cm}
{\centering \begin{tabular}{|c||c|c|c|c|c|c|c|}
\hline 
p\textbackslash{}q&
5&
6&
7&
8&
9&
10&
11\\
\hline 
\hline 
2&
60;12;192&
&
42;6;48&
&
36;4;24&
&
33;3;16\\
\hline 
3&
40;2;16&
&
168;6;128&
32;1;4&
&
120;3;32&
88;2;16\\
\hline 
4&
240;3;64&
&
336;3;128&
&
144;1;8&
&
528;3;128\\
\hline 
5&
&
120;1;4&
840;6;256&
480;3;32&
360;2;32&
&
1320;6;256\\
\hline 
6&
&
&
168;1;8&
&
&
&
264;1;8\\
\hline 
7&
&
&
&
672;3;64&
504;2;32&
840;3;64&
1848;6;256\\
\hline 
8&
&
&
&
&
288;1;4&
&
1056;3;64\\
\hline 
9&
&
&
&
&
&
360;1;4&
792;2;32\\
\hline 
10&
&
&
&
&
&
&
1320;3;64\\
\hline 
\end{tabular}\par}
\vspace{0.3cm}

As an aside, we note that for all minimal models \( \mathcal{M}\left( q,p\right)  \)
\[
N_{0}=\left\{ \begin{array}{cc}
q & p=2\\
4q & p=3\\
4pq & p>3
\end{array}\right. \]
 and

\[
N/N_{0}=\frac{6}{Gcd(6,pq)}\]
for \( p>3 \), while \[
N/N_{0}=\left\{ \begin{array}{cc}
6 & \M{q}{1}{6}\\
3 & \M{q}{4}{6}\\
2 & \M{q}{5}{6}\\
1 & \M{q}{2}{6}
\end{array}\right. \]
for \( p=3 \) and \[
N/N_{0}=\left\{ \begin{array}{cc}
12 & \, \, \, \, \, \, \, \, \, \, \, \, \, \, \, \, \, \M{\, \, \, \, \, \, \, \, q}{1,5,13,17}{24}\\
6 & \, \, \, \, \, \, \, \, \, \, \M{q}{7,23}{24}\\
4 & \, \, \, \, \, \, \, \, \, \M{q}{9,21}{24}\\
3 & \, \, \, \, \, \, \, \, \, \, \, \, \M{q}{11,19}{24}\\
2 & \, \, \M{\, q}{15}{24}\\
1 & \M{q}{3}{24}
\end{array}\right. \]
 for \( p=2 \).

\section{Summary}

As we have seen in the previous sections, the theory of the Galois action supplemented
by the Orbifold Covariance Principle leads to a host of interesting results
about the arithmetic and group theoretic properties of the modular representation.
In particular, we have been able to show that the kernel of the modular representation
is always a congruence subgroup, whose level equals the conductor, i.e. the
order of the Dehn-twist \( T \). Eq.(\ref{kern1}) gives a simple characterization
of the kernel, which allows for efficient algorithms to determine the kernel
explicitly. We have also determined the Galois action on arbitrary modular matrices,
and gave an effective formula to compute the representation matrix of any modular
transformation, Eq.(\ref{modrep1}). Moreover, the results of section 4 put
severe restrictions on the spectrum of the Dehn-twist \( T \), and imply that,
for a given number of primary fields, there is only a finite number of allowed
modular representations. It seems unlikely that the above results could be derived
without using the theory of the Galois action, and provide a nice example of
the use of arithmetic in mathematical physics.

Of course, a host of questions remain open. An obvious one is to find a direct
characterization of the kernel. While Eq.(\ref{kern1}) gives a simple criterion
to decide whether a modular transformation belongs to the kernel in case its
diagonal elements are coprime to the conductor, and in principle the general
case can be reduced to this one, a simple criterion valid for any modular transformation
would be welcome. One might also wonder to find a small set of numerical data
that would characterize the kernel completely. As we have seen, these would
include the parameters \( N \) and \( N_{0} \), which already determine the
kernel to a great extent, but are not enough to characterize it completely.

Another line of study is to investigate the algebro-geometric properties of
the modular curve \( X\left( \mathcal{K}\right)  \) associated to the kernel.
This is a compact Riemann-surface, whose complex structure can be described
by standard techniques, e.g. one can determine the period matrix, etc. It is
also known that the field of meromorphic functions on \( X\left( \mathcal{K}\right)  \)
is a Galois extension of \( \mathbb {C}\left( j\right)  \) with Galois group
\( \Gamma (1)/\mathcal{K} \), where \( j \) is the classical modular invariant.
As the genus one characters of the primaries are meromorphic functions on this
modular curve, it is tempting to expect that some deeper properties of the RCFT
should be encoded in structure of the curve \( X\left( \mathcal{K}\right)  \). 

Still another interesting question concerns the representation \( D \) of \( SL_{2}\left( N\right)  \)
introduced in Section 5. We know already some of its properties, e.g. its kernel
and its behavior under Galois transformations, Eq.(\ref{gal2}), but there is
clearly much more to learn about it. As the representation theory of \( SL_{2}(N) \)
is known, one would be interested to know which irreducible representations
appear in \( D \) and with which multiplicities, etc. Of course, one should
bear in mind that \( D \) is not simply a representation, but a representation
with a choice of a distinguished basis.

Finally, let's make some comments about the relevance of the above on the classification
program. According to Proposition \ref{condbound}, for a given number of primary
fields it is in principle possible to enumerate all matrices \( S \) and \( T \)
that are compatible with the results of this work, i.e. those which can occur
as the generators of the modular representation of a consistent Rational Conformal
Field Theory. Of course, it may happen that such a list would contain spurious
entries which do not correspond to any RCFT, but one may hope that using all
known properties of the modular representation, these spurious solutions can
be discarded, leading to a complete list.

\section{Appendix}

In this appendix we review the definition and the basic properties of the \( \Lambda  \)-matrices
used in Section 3. See \cite{PO2} for more details.

Let \( r=\frac{k}{n} \) be a rational number in reduced form, i.e. with \( n>0 \)
and \( k \) and \( n \) coprime. Choose integers \( x \) and \( y \) such
that \( kx-ny=1 \), and define \( r^{*}=\frac{x}{n} \). Then \( m=\left( \begin{array}{cc}
k & y\\
n & x
\end{array}\right)  \) belongs to \( \Gamma (1) \), and we define the matrix \( \Lambda \left( r\right)  \)
via \[
\Lambda \left( r\right) _{p}^{q}=\omega _{p}^{-r}M_{p}^{q}\omega _{q}^{-r^{*}}\]
or symbolically \[
\Lambda (r)=T^{-r}MT^{-r^{*}}\]
Of course, one should fix some definite branch of the logarithm to make the
above definition meaningful, but different choices lead to equivalent results.

It is a simple matter to show that \( \Lambda (r) \) is well defined, i.e.
does not depend on the actual choice of \( x \) and \( y \). In the same vein,
\( \Lambda (r) \) is periodic in \( r \) with period 1, i.e. \[
\Lambda (r+1)=\Lambda (r)\]

For \( r=0 \) we just get back the \( S \) matrix \[
\Lambda (0)=S\]
and for a positive integer \( n \) we have \begin{equation}
\label{Lambda1n}
\Lambda \left( \frac{1}{n}\right) =T^{-\frac{1}{n}}S^{-1}T^{-n}ST^{-\frac{1}{n}}
\end{equation}

Finally, one can show that \[
\Lambda \left( r^{*}\right) _{p}^{q}=\Lambda (r)_{q}^{p}\]
and \[
\Lambda \left( -r\right) _{p}^{q}=\overline{\Lambda \left( r\right) _{\overline{p}}^{q}}\]

\begin{acknowledgement*}
Many thanks to Terry Gannon for explaining me his work on the Galois action
and the congruence subgroup problem. 
\end{acknowledgement*}


\begin{thebibliography}{10}
\bibitem{Cardy}J. Cardy, Nucl. Phys. \textbf{B270}, 186 (1986).
\bibitem{Ver}E. Verlinde, Nucl. Phys. \textbf{B300}, 360 (1988).
\bibitem{MS}G. Moore and N. Seiberg, Commun. Math. Phys. \textbf{123}, 177 (1989).
\bibitem{TI}P. Bantay, hep-th/0007164.
\bibitem{TFT}E. Witten, Commun. Math. Phys. \textbf{121}, 351 (1989).
\bibitem{Frochlich}J. Frohlich, Int. J. Mod. Phys. \textbf{A20}, 5321 (1989).
\bibitem{Moore}G. Moore, Nucl. Phys. \textbf{B293}, 139 (1987).
\bibitem{Eholzer}W. Eholzer, Commun. Math. Phys. \textbf{172}, 623 (1995).
\bibitem{Eh-Sko}W. Eholzer and N.-P. Skoruppa, Commun. Math. Phys. \textbf{174}, 117 (1995).
\bibitem{BCIR}M. Bauer, A. Coste, C. Itzykson and P. Ruelle, J. Geom. Phys. \textbf{22}, 134
(1997).
\bibitem{CG2}A. Coste and T. Gannon, math-QA/9909080.
\bibitem{BG}J. de Boere and J. Goeree, Commun. Math. Phys. \textbf{139}, 267 (1991).
\bibitem{CG1}A. Coste and T. Gannon, Phys. Lett. \textbf{B323}, 316 (1994).
\bibitem{KS}A. Klemm and M. G. Schmidt, Phys. Lett. \textbf{B245}, 53 (1990).
\bibitem{OCP}P. Bantay, Phys. Lett. \textbf{B488}, 207 (2000).
\bibitem{BHS}L. Borisov, M. B. Halpern and C. Schweigert, Int. J. Mod. Phys. \textbf{A13},
125 (1998).
\bibitem{PO1}P. Bantay, Phys. Lett. \textbf{B419}, 175 (1998).
\bibitem{Bauer}M. Bauer, Advanced Series in Mathematical Physics Vol. 24, 152 (1997).
\bibitem{PO2}P. Bantay, hep-th/9910079.\end{thebibliography}
\end{document}